\documentclass[12pt]{article}
\usepackage{mathrsfs}
\usepackage{amssymb} \textwidth 150mm \textheight 220mm \oddsidemargin
15pt \evensidemargin 0pt \topmargin 0cm \headsep 0.3cm

\usepackage{amsmath}
\usepackage{graphicx}
\usepackage{indentfirst}
\usepackage{cite}
\usepackage{float}
\usepackage{dsfont}
\usepackage{CJK}
\usepackage{multirow}
\usepackage{makecell}

\usepackage{subfigure}

\begin{document}
\title{\Large\bf{The cyclicity of period annulus of cubic isochronous Hamiltonian systems}}
\author{{Jihua Yang\thanks{ E-mail addresses: jihua1113@163.com,\ yangjh@mail.bnu.edu.cn.}}
\\ {\small \it School of Mathematical Sciences, Tianjin Normal University, Tianjin 300387, P. R. China}\\ }
\date{}
\maketitle \baselineskip=0.9\normalbaselineskip \vspace{-3pt}
\noindent
{\bf Abstract}\,  Cima, Ma\~{n}osas and Villadelprat (J. Differ. Equations, 157, 373--413, 1999) proved that a cubic Hamiltonian system possesses an isochronous center at the origin if and only if its Hamiltonian function can be expressed as
\begin{eqnarray*}H_1(x,y)=k_1^2x^2+(k_2y+k_3x+k_4x^2)^2,
\end{eqnarray*}
where $k_1,k_2,k_3,k_4\in\mathbb{R}$, $k_1k_2\neq0$. This paper is devoted to investigating the weak Hilbert's 16th problem  for the
dynamical system associated with the above Hamiltonian function. We show that the maximum number of limit cycles is $n-1$. Furthermore, this number is reached. That is, we solve the weak Hilbert's 16th
problem restricted to cubic Hamiltonian systems with an isochronous center at the origin.
\vskip 0.2 true cm
\noindent

\vskip 0.2 true cm
\noindent

 \section{Introduction and main result}
 \setcounter{equation}{0}
\renewcommand\theequation{1.\arabic{equation}}

For a planar autonomous differential system
\begin{eqnarray}
\frac{dx}{dt}=F_n(x,y),\ \frac{dy}{dt}=G_n(x,y),
\end{eqnarray}
where $F_n(x,y)$ and $G_n(x,y)$ are real polynomials of degree $n$. Let $\Gamma$ be a closed orbit of system (1.1). If there exists a neighborhood of $\Gamma$ containing no other closed orbits of (1.1), then $\Gamma$ is called a limit cycle of (1.1). Geometrically, such a limit cycle is an isolated closed orbit with the distinctive property that all neighboring trajectories asymptotically approach it (either as $t\rightarrow+\infty$ or $t\rightarrow-\infty$) in a spiraling manner, hence the terminology ``limit cycle''. The bifurcation theory of limit cycles in differential systems not only arises naturally in mathematical modeling of biological systems \cite{LMZ,WWW,CYS}, physics \cite{CZ,GH}, economics \cite{P}, mechanics \cite{CDTX}, astronomy  \cite{C}, electronics \cite{BV},
but is also fundamentally connected to Hilbert's 16th problem and its weak version \cite{A,H}.

 At the Second International Congress of Mathematicians held in Paris in 1900, the renowned mathematician D. Hilbert, with his vast knowledge and profound vision, proposed 23 mathematical problems (published in 1902 \cite{H}), among which the second part of the 16th problem asks: what is the least upper bound on the number of limit cycles in system (1.1), and what are their possible relative configurations? In recent decades, mathematicians have conducted extensive research on this problem, yet substantial progress remains limited \cite{I,LJ,HM,LC,HL}. In 1977, Arnold \cite{A} first proposed the weak Hilbert's 16th problem concerning the near-Hamiltonian system:
\begin{equation}
\begin{cases}
\dfrac{dx}{dt} = \dfrac{\partial H(x,y)}{\partial y} + \varepsilon f(x,y), \\[10pt]
\dfrac{dy}{dt} = -\dfrac{\partial H(x,y)}{\partial x} + \varepsilon g(x,y),
\end{cases}
\end{equation}
where $0 < |\varepsilon| \ll 1$, $H(x,y)$ is a polynomial of degree $m+1$, and $f(x,y)$ and $g(x,y)$ are polynomials of degree $n$.
Assume the unperturbed system $(1.2)_{\varepsilon=0}$ possesses a family of closed orbits $\{\Gamma_h\}$, and let $\Sigma$ denote the maximal open interval of $h$ where they exist, i.e.,
$$\Gamma_h = \{(x,y) \mid H(x,y) = h, h \in \Sigma\}.$$
Given that $\Gamma_h$ varies monotonically with $h$, we consider the Abelian integral:
\begin{equation}
I(h) = \oint_{\Gamma_h} g(x,y)dx - f(x,y)dy.
\end{equation}
 The fundamental question is: what is the maximum number  of isolated zeros (counting multiplicities) of the Abelian integral $I(h)$?
 Numerous excellent works have addressed the weak Hilbert's 16th problem; see, for example, \cite{BNY,G,HI98,LZ,IY,II,K,V} and the references therein.

It is well known that the number of limit cycles for a perturbation of a Hamiltonian system is closely related to the number of isolated zeros of the corresponding Abelian integral \cite{G,IY,PL}.  More specifically, the total number of zeros of the Abelian integral $I(h)$, counted with multiplicities, provides an upper bound for the number of limit cycles of system (1.2) bifurcating from the corresponding open period annulus $\bigcup\limits_{h\in\Sigma}\Gamma_h$ \cite{G}. The same is true for the closed period annulus, provided that it is bounded by a homoclinic loop as proved by Roussarie \cite{R}. Moreover, a lower bound for the number of limit cycles is given by the number of multiple simple zeros of $I(h)$.

When $m=n=2$, the weak Hilbert's 16th problem has been completely resolved. This outcome constitutes one of the exceedingly few complete solutions in this research domain, accomplished through more than a decade of sustained scholarly efforts, see \cite{G,HI,LZ,M,ZL}. However, for the case $m=n=3$, only partial results have been obtained so far. Li, Liu and Yang \cite{LLY} proved that there exist polynomials $f(x,y)$ and $g(x,y)$ of degree 3 such that system (1.2) has 13 limit cycles.
Liu and Li \cite{LL}, as well as Yang, Han, Li and Yu \cite{YHLY}, have also established examples demonstrating the existence of 13 limit cycles in cubic differential systems.
For the elliptic Hamiltonian of degree 4 as follows:
\begin{eqnarray*}
H(x,y)=\frac{1}{2}y^2+\frac a4x^4+\frac b3x^3+\frac c2x^2,\ a,b,c\in\mathbb{R}, a\neq0,
\end{eqnarray*}
there are five types of continuous families of ovals on the level curves of $H(x,y)$, depending on the values of the parameters $(a,b,c)$, called the truncated pendulum, the saddle loop, the global center, the cuspidal loop and the figure-eight loop, respectively. When the perturbation is Li\'{e}nard type: $(\alpha+\beta x+\gamma x^2)ydx$, there is a series of papers dealing with the exact number of zeros of the Abelian integrals over five types of ovals.
Horozov \cite{Ho} considered the truncated pendulum, while the seminal work of Dumortier and Li \cite{DL1,DL2,DL3,DL4} addressed fundamental scenarios such as the saddle loop, global center, cuspidal loop, and figure-eight loop, which represent classical results in this domain. Regarding perturbations with $n$-th degree polynomials, Zhao and Zhang \cite{ZZ} proved that the upper bound is $7n+5$. Liu \cite{L} studied the total number of zeros for the ovals in the two annuli surrounded by the figure-eight loop, and improves the upper bound given in \cite{ZZ}. When the Hamiltonian function $H(x, y)$ contains $x^iy^j$, where $i$ and $j$ are positive
integers, Zhou and Li \cite{ZL1} obtained the algebraic structure of the Abelian integral for the Hamiltonian
$$H(x, y) = x^2 + y^2 + ax^4 + bx^2y^2 + cy^4,\ a,b,c\in\mathbb{R},$$
and an upper bound of the number of zeros of Abelian integral $I(h)$ was given for a special
case $a>0, b=0$ and $c=1$. Later, Chen and Yu \cite{CY} obtained an upper bound for $a,b,c\in\mathbb{R}$. Wu, Zhang and Li \cite{WZL} obtained an upper bound for the case
$$H(x, y) = x^2 + y^2 + cx^2y^2-x^4 + y^4,\ c > -2.$$
Yang and Zhao \cite{YZ} gave an upper bound (except the butterfly phase portrait) for the case
$$H(x, y) = -x^2 + ax^4+ bx^2y^2 + cy^4,\ a,b,c\in\mathbb{R}, c \neq 0.$$
Later, Yang, Sui and Zhao \cite{YSZ} got an upper bound of the above system with the butterfly phase portrait.
Chang, Zhao and Wang \cite{CZW} derived an upper bound  for the case
$$H(x, y) = \alpha x^2 +\beta y^2 +ax^4 + bx^2y^2+ cy^4,\ \alpha,\beta, a,b,c\in\mathbb{R}, \alpha\beta<0.$$

It is worth noting that the systems studied in the aforementioned literature are all symmetric with respect to the $x$-axis or $y$-axis, which reduces the number of generators for the Abelian integrals. Inspired by these works, this paper focuses on the limit cycle bifurcations in a class of cubic Hamiltonian systems that lack symmetry about the coordinate axes. In 1999, Cima, Ma\~{n}osas and Villadelprat \cite{CMV} determined all the cubic Hamiltonian systems that have an isochronus center at the origin. They proved the following conclusion:
\vskip 0.2 true cm

\noindent
{\bf Theorem 1.1}\,\cite{CMV} {\it A cubic Hamiltonian system has an isochronous center at the origin if and only if after a linear change of coordinates its Hamiltonian function can be written as
\begin{eqnarray}
H_1(x,y)=k_1^2x^2+(k_2y+k_3x+k_4x^2)^2,
\end{eqnarray}
where $k_i\in\mathds{R}$ for $i=1,2,3,4$ and $k_1$ and $k_2$ are different from zero.}
\vskip 0.2 true cm

In the present paper, we study the weak Hilbert's 16th problem for the dynamical system associated with Hamiltonian function (1.4). The Hamiltonian system corresponding to (1.4) is
\begin{eqnarray}\begin{cases}
\frac{dx}{dt}=-2k_2(k_2y+k_3x+k_4x^2),\\
\frac{dy}{dt}=2k_1^2x+2(k_3+2k_4x)(k_2y+k_3x+k_4x^2),
\end{cases}\end{eqnarray}
Letting $t_1=2k_2k_3t$, $x_1=\frac{k_4}{k_3}x$ and $y_1=\frac{k_2k_4}{k_1^2+k_3^2}y$ $(k_i\neq0,i=1,2,3,4)$, one can change system (1.5) into
\begin{eqnarray}\begin{cases}
\frac{dx}{dt}=-\lambda^{-1}y-x-x^2,\\
\frac{dy}{dt}=x+y+2xy+3\lambda x^2+2\lambda x^3,
\end{cases}\end{eqnarray}
with the Hamiltonian function
\begin{eqnarray}
H(x,y)=\frac{1}{2}x^2+\lambda x^3+\frac{1}{2}\lambda x^4+\frac{1}{2}\lambda^{-1}y^2+xy+x^2y,
\end{eqnarray}
where $\lambda=\frac{k_3^2}{k_1^2+k_3^2}$. Clearly, $0<\lambda<1$. Here and below, we shall omit the subscript 1. System (1.6) has an isochronous center at the origin and  a family of periodic orbits, denoted by $$\Gamma_h=\{(x,y):H(x,y)=h,h\in(0,+\infty)\}.$$
The parabola $y=-\lambda x^2-\lambda x$ divides $\Gamma_h$ into an upper arc and a lower arc, with their respective function expressions given by
$$y=-\lambda x^2-\lambda x+\sqrt{(\lambda^2-\lambda)x^2+2\lambda h}$$ and $$y=-\lambda x^2-\lambda x-\sqrt{(\lambda^2-\lambda)x^2+2\lambda h}.$$
The coordinates of the two intersection points between the parabola (the red curve in Fig. 1) and $\Gamma_h$ (the blue closed curves in Fig. 1) are as follows:
$$\Big(\sqrt{\frac{2h}{1-\lambda}},-\lambda \sqrt{\frac{2h}{1-\lambda}}-\frac{2\lambda h}{1-\lambda}\Big),\ \Big(-\sqrt{\frac{2h}{1-\lambda}},\lambda \sqrt{\frac{2h}{1-\lambda}}-\frac{2\lambda h}{1-\lambda}\Big).$$
Our main result is the following theorem.
\vskip 0.2 true cm

\noindent
{\bf Theorem 1.2}\, {\it Consider the following perturbation of system (1.6):\begin{eqnarray}
\begin{cases}
\frac{dx}{dt}=-\lambda^{-1}y-x-x^2+\varepsilon\sum\limits_{i+j=0}^na_{i,j}x^iy^j,\\
\frac{dy}{dt}=x+y+2xy+3\lambda x^2+2\lambda x^3+\varepsilon\sum\limits_{i+j=0}^nb_{i,j}x^iy^j,\\
\end{cases}
\end{eqnarray}
where $0<|\varepsilon|\ll1.$
  Then, by using the Abelian integral, the upper bound for the number of limit cycles of system (1.8) bifurcating from the period annulus is $n-1$ for $n\geq2$, counted with multiplicities. Moreover, this bound is sharp.}\vskip 0.2 true cm

\vskip 0.2 true cm

\noindent
{\bf Remark 1.1}\, (i) One major challenge in this paper lies in analyzing the algebraic structure of the Abelian integral $I(h)$. As demonstrated, the number of generators of $I(h)$ depends on the degree $n$ of the perturbation polynomials, which constitutes the key distinction from existing literature. To address this difficulty, we classify the terms $I_{i,j}(h)$ appearing in $I(h)$ into two categories:
\vskip 0.2 true cm

(a) formula-iterable terms admitting recursive computation: $I_{i,j}(h),i\geq2,j\geq1$;
\vskip 0.2 true cm

(b) non-iterable terms requiring alternative treatment: $I_{i,j}(h),i=0,1,j\geq1$.\vskip 0.2 true cm

\noindent
For the non-iterable terms $I_{i,j}(h)$, we first derive the differential equations they satisfy and then obtain their explicit expressions by solving these differential equations.
\vskip 0.2 true cm

\noindent
(ii) After obtaining the explicit expression of $I(h)$, verifying the linear independence of its coefficients becomes essential for determining the lower bound of the number of limit cycles. This constitutes another fundamental challenge in our work, which we successfully overcome through an innovative application of mathematical induction.
\vskip 0.2 true cm

\noindent
(iii) As shown in Fig. 1, the phase portrait of system (1.6) exhibits no symmetry whatsoever--neither about the coordinate axes nor about the origin. This inherent asymmetry inevitably leads to a larger number of generators for the corresponding Abelian integral $I(h)$  than classical methodologies can accommodate, rather than a restriction to a small finite set (e.g., two or three generators).

\vskip 0.2 true cm

The paper is organized as follows. The detailed expression of the Abelian integral $I(h)$ is obtained in Section 2. The proof of the Theorem 1.1 and some numerical simulations are presented in Section 3. The discussion is then presented in the final section.

 \section{The algebraic structure of Abelian integral}
 \setcounter{equation}{0}
\renewcommand\theequation{2.\arabic{equation}}
For abbreviation we denote
$$I_{i,j}(h)=\oint_{\Gamma_h}x^iy^jdx,\ h\in(0,+\infty),\ i,j\in\mathds{N}.$$
It is straightforward to check that $I_{n,0}(h)=0$. Direct computation by applying Green's formula yields
\begin{eqnarray}\begin{aligned}
I(h)=&\sum\limits_{i+j=0}^nb_{i,j}\oint_{\Gamma_h}x^iy^jdx-\sum\limits_{i+j=0}^na_{i,j}\oint_{\Gamma_h}x^iy^jdy\\
=&\sum\limits_{i+j=1,j\geq1}^n\xi_{i,j}I_{i,j}(h),
\end{aligned}\end{eqnarray}
in view of
\begin{eqnarray}\begin{aligned}
\oint_{\Gamma_h}x^iy^jdy=-\frac{i}{j+1}I_{i-1,j+1}(h),
\end{aligned}\end{eqnarray}
where $\xi_{i,j}=b_{i,j}+\frac{i+j}{j}a_{i+1,j-1}$ and can be chosen as free parameters.
\vskip 0.2 true cm

\noindent
{\bf Lemma 2.1}\, {\it The following relationship holds
\begin{eqnarray}
I_{n,1}(h)=\begin{cases}0,\qquad\qquad\qquad\qquad\qquad\quad\quad\qquad\ \  n\ odd,\\
-4\int_0^{\sqrt{\frac{2h}{1-\lambda}}}x^n\sqrt{(\lambda^2-\lambda)x^2+2\lambda h}dx,\ n\ even.
\end{cases}
\end{eqnarray} }
\vskip 0.2 true cm

\noindent
{\bf Proof}\,  Some direct computation yields
$$\begin{aligned}I_{n,1}(h)=&\oint_{\Gamma_h}x^nydx=\int_{\sqrt{\frac{2h}{1-\lambda}}}^{-\sqrt{\frac{2h}{1-\lambda}}}x^n\big[-\lambda x^2-\lambda x+\sqrt{(\lambda^2-\lambda)x^2+2\lambda h}\big]dx\\&+
\int_{-\sqrt{\frac{2h}{1-\lambda}}}^{\sqrt{\frac{2h}{1-\lambda}}}x^n\big[-\lambda x^2-\lambda x-\sqrt{(\lambda^2-\lambda)x^2+2\lambda h}\big]dx\\
=&-2\int_{-\sqrt{\frac{2h}{1-\lambda}}}^{\sqrt{\frac{2h}{1-\lambda}}}x^n\sqrt{(\lambda^2-\lambda)x^2+2\lambda h}dx.
\end{aligned}$$
Note that when $n$ is odd, the integrand of the above integral is an odd function, and when $n$ is even, the integrand is an even function. Then (2.3) follows immediately by symmetry. This completes the proof.\quad $\lozenge$\vskip 0.2 true cm

The following lemma plays a crucial role in determining  the algebraic structure of the Abelian integral $I(h)$.

\vskip 0.2 true cm

\noindent
{\bf Lemma 2.2}\, {\it For $n\geq5$, the Abelian integral $I(h)$ can be expressed as
\begin{eqnarray}\begin{aligned}
I(h)=&\sum\limits_{i=1}^{n-2}\bar{P}_{[\frac{n+1-i}{4}]}(h)I_{0,i}(h)+\alpha_1I_{0,n-1}(h)+\alpha_2I_{0,n}(h)\\&+
\sum\limits_{i=2}^{n-3}\bar{Q}_{[\frac{n+1-i}{4}]}(h)I_{1,i}(h)+\beta_1I_{1,n-2}(h)+\beta_2I_{1,n-1}(h),
\end{aligned}\end{eqnarray}
where $\bar{P}_l(h)$ and $\bar{Q}_l(h)$ are polynomials of degree $l$. }
\vskip 0.2 true cm

\noindent
{\bf Proof}\, Differentiating both sides of $H(x,y)=h$  in (1.7) with respect to $x$ gives
\begin{eqnarray}\begin{aligned}
x+y+2xy+3\lambda x^2+2\lambda x^3+\lambda^{-1}y\frac{\partial y}{\partial x}+x\frac{\partial y}{\partial x}+x^2\frac{\partial y}{\partial x}=0.
\end{aligned}\end{eqnarray}
Multiplying both sides of (2.5) by $x^{i-3}y^jdx$ and integrating along $\Gamma_h$, one gets
\begin{eqnarray}\begin{aligned}
I_{i,j}(h)=&\frac{1}{2\lambda}\Big[\frac{i-3}{\lambda(j+2)}I_{i-4,j+2}(h)-I_{i-2,j}(h)+\frac{i-2j-3}{j+1}I_{i-2,j+1}(h)\\&-3\lambda I_{i-1,j}(h)+\frac{i-j-3}{j+1}I_{i-3,j+1}(h)\Big].
\end{aligned}\end{eqnarray}
At the same time, multiplying both sides of $H(x,y)=h$ by $x^iy^{j-2}dx$ and integrating along $\Gamma_h$ imply
\begin{eqnarray}\begin{aligned}
I_{i,j}(h)=&2\lambda hI_{i,j-2}(h)-\lambda I_{i+2,j-2}(h)-2\lambda I_{i+2,j-1}(h)\\&-\lambda^2I_{i+4,j-2}(h)-2\lambda I_{i+1,j-1}(h)-2\lambda^2 I_{i+3,j-2}(h).
\end{aligned}\end{eqnarray}
On account of (2.6) and (2.7), one can derive two recurrence formulas
\begin{eqnarray}\begin{aligned}
I_{i,j}(h)=&\frac{1}{\lambda(i+2j+1)}\Big[2(i-3)hI_{i-4,j}(h)-(i+j-1)I_{i-2,j}(h)\\&-\lambda(2i+3j)I_{i-1,j}(h)-\frac{j(i+2j+1)}{j+1}I_{i-2,j+1}(h)\\&-\frac{j(i+j-1)}{j+1}I_{i-3,j+1}(h)\Big],
\end{aligned}\end{eqnarray}
and
\begin{eqnarray}\begin{aligned}
I_{i,j}(h)=&\frac{\lambda j}{i+2j+1}\Big[4hI_{i,j-2}(h)- I_{i+2,j-2}(h)-\lambda I_{i+3,j-2}(h)\\&-\frac{i+2j+1}{j-1} I_{i+2,j-1}(h)-\frac{i+3j-1}{j-1} I_{i+1,j-1}(h)\Big].
\end{aligned}\end{eqnarray}
Eliminating $I_{i+2,j-1}(h)$ and $I_{i+3,j-2}(h)$ in (2.9) using (2.8) leads to
\begin{eqnarray}\begin{aligned}
&4\lambda j(j-1)(i+2j)hI_{i,j-2}(h)+j(i+j)(i+2j+\lambda j-2\lambda)I_{i,j-1}(h)\\&-2j(i-1)(i+2j)hI_{i-2,j-1}(h)+(i+j)(j-1)(i+2j)I_{i-1,j}(h)\\&
+\lambda j(i+j)(i+2j)I_{i+1,j-1}(h)-\lambda j(j-1)(i+2j-2\lambda i-3\lambda j)I_{i+2,j-2}(h)\\&
-2\lambda ij(j-1)hI_{i-1,j-2}(h) +\lambda j(j-1)(i+j)I_{i+1,j-2}(h)=0.
\end{aligned}\end{eqnarray}
Taking $i\rightarrow i-1,j\rightarrow j+1$ in (2.10) yields another required recurrence formula
\begin{eqnarray}\begin{aligned}
I_{i,j}(h)=&-\frac{1}{\lambda (j+1)(i+j)(i +2 j +1)}\big[4\lambda j(j +1)(i+2j+1) h I_{i-1, j -1}(h) \\
&-2(i -2)(j+1)(i+2j+1)h I_{i-3,j}(h)+j(i+j)(i+2j+1)I_{i-2,j+1}(h)\\&
+(j+1)(i+j)\big(i+2j+1+\lambda j-\lambda\big)I_{i-1,j}(h)\\&
-\lambda j(j+1)\big(i+2j+1-2\lambda i-3\lambda j-\lambda\big)I_{i+1,j-1}(h)\\&
-2\lambda j(i-1)(j+1)hI_{i-2,j-1}(h)+\lambda j(j+1)(i+j)I_{i,j-1}(h)
\big].
\end{aligned}\end{eqnarray}

We are now in a position to prove (2.4) by induction on $n$ using the recurrence formulas (2.8) and (2.11). It follows from (2.8), (2.11) and Lemma 2.1 that
\begin{eqnarray}\begin{aligned}
I_{2,1}(h)=&-\frac{1}{2\lambda}I_{0,2}(h),\\
I_{2,2}(h)=&\frac{1}{7}hI_{0,1}(h)-\frac{\lambda+6}{7\lambda}I_{0,2}(h)-\frac{2}{3\lambda}I_{0,3}(h),\\
I_{3,2}(h)=&-\frac{3}{14}hI_{0,1}(h)-\frac{3}{14\lambda}I_{0,2}(h)+\frac{2}{3\lambda}I_{0,3}(h)+\frac{3\lambda+14}{14\lambda}I_{1,2}(h)-\frac{2}{3\lambda}I_{1,3}(h),\\
I_{4,1}(h)=&-\frac{1}{28\lambda}h I_{0,1}(h)+\frac{3}{14\lambda^2}I_{0,2}(h)+\frac{1}{3\lambda^2}I_{0,3}(h)+\frac{\lambda+7}{28\lambda^2}I_{1,2}(h),\\
I_{2,3}(h)=&\frac{42\lambda-37}{210}hI_{0,1}(h)+\frac{1}{210\lambda}(28\lambda h+42\lambda-37)I_{0,2}(h)\\&-\frac{28}{45\lambda}(\lambda-1)I_{0,3}(h)-\frac{3}{4\lambda}I_{0,4}(h)+\frac{1}{5\lambda}(2\lambda-7)I_{1,3}(h)
\\&-\frac{1}{210\lambda}(504\lambda h+42\lambda^2+159\lambda-196)I_{1,2}(h),
\end{aligned}\end{eqnarray}
which yields that (2.4) is valid for $n=5$.  Some tedious manipulation using (2.8) and (2.11) gives rise to
\begin{eqnarray}
\mathbf{\Phi}\left(\begin{matrix}
                 I_{2,n-2}(h)\\
                 I_{3,n-3}(h)\\
                  I_{4,n-4}(h)\\
                 \vdots\\
                  I_{n-2,2}(h)\\
                  I_{n-1,1}(h)
                \end{matrix}\right)=\left(\begin{matrix}
                 v_{1}(h)\\
                 v_{2}(h)\\
                  v_{3}(h)\\
                 \vdots\\
                  v_{n-3}(h)\\
                  v_{n-2}(h)
                \end{matrix}\right),
\end{eqnarray}
where
\begin{eqnarray*}\begin{aligned}
v_1(h)=&-\frac{1}{\lambda n(n-1)(2n-1)}\big[4\lambda(n - 1)(n - 2)(2n - 1)hI_{1, n - 3}(h)\\&+n(n - 2)(2n - 1)I_{0, n - 1}(h) + n(n-1)((n - 3)\lambda - 1 + 2n )I_{1, n - 2}(h)\\&
- 2\lambda(n - 1)(n - 2)hI_{0, n - 3}(h)+\lambda n(n - 1)(n - 2)I_{2, n - 3}(h)\big],  \\
v_2(h)=&-\frac{1}{2\lambda(n-2)}\big[(n-2)I_{1,n-3}(h)+2(n-3)I_{1,n-2}(h) \\&+3\lambda(n-2)I_{2,n-3}(h) +(n-3)I_{0,n-2}(h) \big],\\
v_3(h)=&-\frac{1}{\lambda(n-3)(2n-3)}\big[(n^2-5n+4)I_{1,n-3}(h)+(n^2-4n+3)I_{2,n-4}(h) \\&+(2n^2-11n+12)I_{2,n-3}(h) -2(n-3)hI_{0,n-4}(h)\\&+\lambda(n-3)(3n-4)I_{3,n-4}(h)  \big],\\
&\qquad\qquad\qquad\qquad\qquad\qquad\vdots\\
v_{n-3}(h)=&\frac{1}{3\lambda (n+3)}\big[6(n - 5)hI_{n - 6, 2}(h) -2 (n -1)I_{n - 5, 3}(h) -3 (n -1)I_{n - 4, 2}(h)\\& -2(n +3)I_{n - 4, 3}(h) - 6\lambda(n + 1)I_{n - 3, 2}(h)\big],\\
v_{n-2}(h)=&\frac{1}{2\lambda (n+2)}\big[4(n - 4)hI_{n - 5, 1}(h) - (n -1)I_{n - 4, 2}(h) -2(n -1)I_{n - 3, 1}(h)\\& -(n +2)I_{n - 3, 2}(h) - 2\lambda(2n+1)I_{n - 2, 1}(h)\big],
\end{aligned}\end{eqnarray*}
$$\mathbf{\Phi}=\left(\begin{matrix}
                 1&\frac{(n-2)((3n - 1)\lambda + 1 - 2n)}{n(2n-1)}&0&\cdots&0&0\\
                 0&1&0&\cdots&0&0\\
                  0&0&1&\cdots&0&0\\
                \vdots&\vdots&\vdots&\ddots&\vdots&\vdots\\
                0&0&0&\cdots&1&0\\
                  0&0&0&\cdots&0&1\\
                \end{matrix}\right).$$
 In order to invoke the induction hypothesis, we rewrite $I(h)$ in the following form:
\begin{eqnarray*}\begin{aligned}
I(h)=\sum\limits_{i+j=0}^n\xi_{i,j}I_{i,j}(h)=&\sum\limits_{i+j=0}^{n-1}\xi_{i,j}I_{i,j}(h)+\xi_{0,n}I_{0,n}(h)+\xi_{1,n-1}I_{1,n-1}(h)\\&
+\xi_{2,n-2}I_{2,n-2}(h)+\cdots+\xi_{n-1,1}I_{n-1,1}(h).
\end{aligned}\end{eqnarray*}
Substituting (2.13) into the above equation and applying the induction hypothesis immediately yields (2.4). This completes the proof.\quad $\lozenge$
\vskip 0.2 true cm

\noindent
{\bf Remark 2.1} (i) Although the recurrence formulas (2.8) and (2.11) play a crucial role in studying the algebraic structure of the Abelian integral
$I(h)$, $I_{0,n}(h)$ and $I_{1,n-1}(h)$ cannot be iterated using these two formulas. This implies that $I(h)$ cannot be represented by a finite set of generators, which constitutes the most significant distinction from previous literature.

\vskip 0.2 true cm

\noindent
(ii) As evident from the proof of Lemma 2.2, both the iterative formula (2.11) itself and its derivation process are remarkably complex, owing to the intricate nature of the first integral of system (1.6).

\vskip 0.2 true cm
In the lemma below, we present the exact expressions of $I(h)$ corresponding to $n=1,2,3,4,$ which are derived by direct computation using Lemma 2.1.
\vskip 0.2 true cm

\noindent
{\bf Lemma 2.3} {\it The Abelian integral $I(h)$ can be written as
\begin{eqnarray*}I(h)=\begin{cases}
\xi_{0,1}I_{0,1}(h),\qquad\qquad\qquad\qquad\qquad \qquad\qquad\qquad\qquad\qquad\  n=1,\\
\xi_{0,1}I_{0,1}(h)+\xi_{0,2}I_{0,2}(h),\qquad\qquad\qquad\qquad\qquad\qquad\qquad\ \ \  n=2,\\
\xi_{0,1}I_{0,1}(h)+(\xi_{0,2}-\frac{1}{2\lambda}\xi_{2,1})I_{0,2}(h)+\xi_{0,3}I_{0,3}(h)+\xi_{1,2}I_{1,2}(h), \ n=3,\\
\xi_{0,1}I_{0,1}(h)+(\xi_{0,2}-\frac{1}{2\lambda}\xi_{2,1}-\frac{1}{2\lambda}\xi_{2,2})I_{0,2}(h)+\xi_{0,3}I_{0,3}(h)\\ \quad\qquad\quad\
+\xi_{1,2}I_{1,2}(h)+\xi_{0,4}I_{0,4}(h)+\xi_{1,3}I_{1,3}(h),\qquad\qquad \ \,\,  n=4.
\end{cases}
\end{eqnarray*}}

Our current objective is to compute the integrals $I_{0,i}(h)$ for $i=1,2,\cdots,n$, along with $I_{1,j}(h)$ for $j=2,3,\cdots,n-1$ in Lemmas 2.2 and 2.3. This will enable us to derive a more detailed expression for the Abelian integral $I(h)$. The following lemma plays a crucial role in the computation of $I_{0,i}(h)$ and $I_{1,j}(h)$.
\vskip 0.2 true cm

\noindent
{\bf Lemma 2.4} {\it For $n\geq3$, the following equalities hold:
\begin{eqnarray}\begin{aligned}
I_{2,n}(h)=&\sum\limits_{i=1}^{n-1}\bar{P}^i_1(h)I_{0,i}(h)+\bar{\alpha}_1I_{0,n}(h)+\bar{\alpha}_2I_{0,n+1}(h)\\&+\sum\limits_{i=2}^{n-1}\bar{Q}^i_1(h)I_{1,i}(h)+\bar{\beta}_1I_{1,n}(h),\\
I_{3,n}(h)=&\sum\limits_{i=1}^{n-1}\tilde{P}^i_1(h)I_{0,i}(h)+\tilde{\alpha}_1I_{0,n}(h)+\tilde{\alpha}_2I_{0,n+1}(h)\\&+\sum\limits_{i=2}^{n-1}\tilde{Q}^i_1(h)I_{1,i}(h)+\tilde{\beta}_1I_{1,n}(h)+\tilde{\beta}_2I_{1,n+1}(h),\\
I_{4,n}(h)=&\sum\limits_{i=1}^{n}\hat{P}^i_1(h)I_{0,i}(h)+\hat{\alpha}_1I_{0,n+1}(h)+\hat{\alpha}_2I_{0,n+2}(h)\\&+\sum\limits_{i=2}^{n}\hat{Q}^i_1(h)I_{1,i}(h)+\hat{\beta}_1I_{1,n+1}(h),\\
\end{aligned}\end{eqnarray}
where $\bar{P}_1^i(h)$, $\bar{Q}_1^i(h)$, $\tilde{P}_1^i(h)$, $\tilde{Q}_1^i(h)$, $\hat{P}_1^i(h)$ and $\hat{Q}_1^i(h)$ are linear polynomials of $h$ and $\bar{\alpha}_i,\bar{\beta}_i,\tilde{\alpha}_i,\tilde{\beta}_i,\hat{\alpha}_i,\hat{\beta}_i\in\mathbb{R}$.}
 \vskip 0.2 true cm

\noindent
{\bf Proof}\, We only prove the first relation in (2.14) by mathematical induction, using (2.8) and (2.11). The proofs of the other two relations proceed in a similar fashion. The fifth relation in (2.12) implies that the conclusion holds when $n=3$. Taking $(i,j)=(2,n+1)$ in (2.11) and $(i,j)=(3,n)$ in (2.8) give rise to
\begin{eqnarray}\begin{aligned}
I_{2,n+1}(h)=&\frac{2(n+1)}{(n+3)(2n+5)}hI_{0,n}(h)-\frac{4(n+1)}{n+3}hI_{1,n}(h)-\frac{n+1}{2n+5}hI_{2,n}(h)\\&-\frac{n+1}{\lambda(n+2)}I_{0,n+2}(h)
-\frac{2n+\lambda n+5}{\lambda(2n+5)}I_{1,n+1}(h)\\&+\frac{(n+1)(2n-3\lambda n-8\lambda+5)}{2n^2+11n+15}I_{3,n}(h),
\end{aligned}\end{eqnarray}
and
\begin{eqnarray}\begin{aligned}
I_{3,n}(h)=&-\frac{1}{2\lambda}I_{1,n}(h)-\frac{3}{2}I_{2,n}(h)-\frac{n}{\lambda(n+1)}I_{1,n+1}(h)-\frac{n}{2\lambda(n+1)}I_{0,n+1}(h).
\end{aligned}\end{eqnarray}
By substituting (2.16) into (2.15), one obtains
\begin{eqnarray}\begin{aligned}
I_{2,n+1}(h)=&\frac{(n+1)(9\lambda n-8n+24\lambda-21)}{4n^2+22n+30}I_{2,n}(h)+\frac{2(n+1)}{(n+3)(2n+5)}hI_{0,n}(h)\\&+\frac{n(3\lambda n-2n+8\lambda-5)}{2\lambda(n+3)(2n+5)}I_{0,n+1}(h)-\frac{n+1}{\lambda(n+2)}I_{0,n+2}(h)\\&-\frac{(n+1)\big((16\lambda n+40\lambda)h-3\lambda n+2n-8\lambda+5\big)}{2\lambda(n+3)(2n+5)}I_{1,n}(h)
\\&+\frac{\lambda n- 2n-3}{\lambda(n+3)}I_{1,n+1}(h).
\end{aligned}\end{eqnarray}
The conclusion is immediately established by (2.17) together with the induction hypothesis. This completes the proof.\quad $\lozenge$
\vskip 0.2 true cm

In order to determine $I_{0,i}(h)$ and $I_{1,j}(h)$, in addition to Lemma 2.4, we also need to find the differential equations they satisfy, as provided by the following lemma.
\vskip 0.2 true cm

\noindent
{\bf Lemma 2.5} {\it For $n\geq2$, the following differential equations hold:
 \begin{eqnarray}\begin{aligned}
&I_{0,n}(h)=4hI'_{0,n}(h)-\frac{n}{n+1}I'_{1,n+1}(h)-I'_{2,n}(h)-\lambda I'_{3,n}(h),\\
\end{aligned}\end{eqnarray}
 \begin{eqnarray}\begin{aligned}
I_{1,n}(h)=&\frac{2(n+3)}{3}hI'_{1,n}(h)-\frac{n+3}{6}I'_{3,n}(h)\\&-\frac{\lambda(n+3)}{6}I'_{4,n}(h)-\frac{n(n+3)}{6n+6}I'_{2,n+1}(h),
\end{aligned}\end{eqnarray}where $'$ means a differentiation with respect to $h$.}
\vskip 0.2 true cm

\noindent
{\bf Proof}\, In the equation $H(x,y)=h$, we regard $y$ as a bivariate function of $x$ and $h$. Differentiating both sides of $H(x,y)=h$ with respect to $h$ gives $$\frac{\partial y}{\partial h}=\frac{1}{x+x^2+\lambda^{-1}y},$$ which yields that
\begin{eqnarray}
I'_{i,j}(h)=j\oint_{\Gamma_h}\frac{x^iy^{j-1}}{x+x^2+\lambda^{-1}y}dx.
\end{eqnarray}
A straightforward computation using (2.20), one has
\begin{eqnarray}
I_{i,j}(h)=\frac{1}{\lambda (j+2)}I'_{i,j+2}(h)+\frac{1}{j+1}I'_{i+1,j+1}(h)+\frac{1}{j+1}I'_{i+2,j+1}(h).
\end{eqnarray}
Multiplying both sides of (2.20) by $h$, one obtains
\begin{eqnarray}\begin{aligned}
hI'_{i,j}(h)=&\frac{1}{2}I'_{i+2,j}(h)+ \lambda I'_{i+3,j}(h)+\frac{\lambda}{2}I'_{i+4,j}(h)+\frac{j}{2\lambda (j+2)}I'_{i,j+2}(h)\\&+\frac{j}{j+1}I'_{i+1,j+1}(h)+\frac{j}{j+1}I'_{i+2,j+1}(h).
\end{aligned}\end{eqnarray}

From another perspective, some routine calculations using (2.2) show
\begin{eqnarray}\begin{aligned}
I_{i,j}(h)=&-\frac{j}{i+1}\oint_{\Gamma_h}x^{i+1}y^{j-1}dy\\
=&\frac{1}{i+1}I'_{i+2,j}(h)+\frac{j}{(i+1)(j+1)}I'_{i+1,j+1}(h)+\frac{3\lambda}{i+1}I'_{i+3,j}(h)
\\&+\frac{2j}{(i+1)(j+1)}I'_{i+2,j+1}(h)+\frac{2\lambda}{i+1}I'_{i+4,j}(h).
\end{aligned}\end{eqnarray}
It follows from (2.21), (2.22) and (2.23) that
{\small\begin{eqnarray}
I_{i,j}(h)=\frac{i^2+3i+2ij+2}{(i+1)^2 (i+2)}\Big(4hI'_{i,j}(h)-I'_{i+2,j}(h)-\lambda I'_{i+3,j}(h)-\frac{j}{j+1}I'_{i+1,j+1}(h)\Big).
\end{eqnarray}}
Taking $(i,j)=(0,n),(1,n)$ in (2.24) gives (2.18) and (2.19). This completes the proof.\quad $\lozenge$

\vskip 0.2 true cm

\noindent
{\bf Lemma 2.6} {\it Let $n$ be a positive integer with $n\geq2$. Then we have
%\begin{eqnarray}\begin{aligned}
%I_{0,n}(h)=\begin{cases}
%h^2P_{n-2}(h),\ \ n\ \textup{even},\\
%hP_{n-1}(h),\quad\ n\ \textup{odd},\\
%\end{cases}
%I_{1,n}(h)=\begin{cases}
%h^2Q_{n-2}(h),\ \ n\ \textup{even},\\
%hQ_{n-1}(h),\quad\ n\ \textup{odd},\\
%\end{cases}\end{aligned}\end{eqnarray}
\begin{eqnarray}
I_{0,n}(h)=hP_{n-1}(h),\ \ I_{1,n}(h)=hQ_{n-1}(h),
\end{eqnarray}
where $P_{n-1}(h)$ and $Q_{n-1}(h)$ are polynomials of degree $n-1$.}
\vskip 0.2 true cm

\noindent
{\bf Proof}\, A direct calculation gives rise to
\begin{eqnarray}\begin{aligned}
%&I_{0,1}(h)=-4\int_0^{\sqrt{\frac{2h}{1-\lambda}}}\sqrt{2\lambda h+(\lambda^2-\lambda)x^2}dx=-\frac{2\pi\sqrt{\lambda}}{\sqrt{1-\lambda}} h,\\
&I_{0,2}(h)=I_{1,2}(h)=8\lambda\int_0^{\sqrt{\frac{2h}{1-\lambda}}}x^2\sqrt{2\lambda h+(\lambda^2-\lambda)x^2}dx=\frac{2\pi\lambda^\frac32}{(1-\lambda)^\frac32}  h^2,\\
&I_{0,3}(h)=-4\lambda\int_0^{\sqrt{\frac{2h}{1-\lambda}}}(3\lambda x^4+4\lambda x^2-x^2+2h)\sqrt{2\lambda h+(\lambda^2-\lambda)x^2}dx\\&
\qquad\ \ =-\frac{3\pi\lambda^\frac32}{(1-\lambda)^\frac52}h(\lambda h^2-\lambda h+h),\\
&I_{1,3}(h)=-24\lambda^2\lambda\int_0^{\sqrt{\frac{2h}{1-\lambda}}}x^4\sqrt{2\lambda h+(\lambda^2-\lambda)x^2}dx=-\frac{6\pi\lambda^\frac52}{(1-\lambda)^\frac52}  h^3,\
\end{aligned}\end{eqnarray}
which imply that (2.25) is valid for $n=2,3$. Assume that (2.25) holds for all $m\leq n$. In order to find $I_{0,n+1}(h)$ and $I_{1,n+1}(h)$, we proceed to establish the differential equations they satisfy.
 In fact, replacing $n$ by $n-1$ in (2.19), one has
 \begin{eqnarray}\begin{aligned}
I_{1,n-1}(h)=&\frac{2(n+2)}{3}hI'_{1,n-1}(h)-\frac{n+2}{6}I'_{3,n-1}(h)\\&-\frac{\lambda(n+2)}{6}I'_{4,n-1}(h)-\frac{(n-1)(n+2)}{6n}I'_{2,n}(h).
\end{aligned}\end{eqnarray}
 Plugging (2.14) into $(2.27)$ and applying the induction hypothesis implies a simple differential equation satisfied by $I'_{0,n+1}(h)$
 \begin{eqnarray}
 I'_{0,n+1}(h)=\check{P}_{n}(h),
 \end{eqnarray}
 where $\check{P}_{n}(h)$ is a polynomial of degree $n$. Solving differential equation (2.28) yields the first equality in (2.25). Substituting (2.14) into (2.18) and applying the induction hypothesis together with (2.28), one can obtain the differential equation satisfied by $I'_{1,n+1}(h)$ as follows
\begin{eqnarray}
 I'_{1,n+1}(h)=\check{Q}_{n}(h),
 \end{eqnarray}
 where $\check{Q}_{n}(h)$ is a polynomial of degree $n$.  Solving the above differential equation gives the second equality in (2.25). This completes the proof.\quad $\lozenge$
\vskip 0.2 true cm

With the preceding preparations in place, we can now derive a more complete expression for the Abelian integral $I(h)$, which in fact takes polynomial form.
\vskip 0.2 true cm

\noindent
{\bf Proposition 2.1} {\it The  Abelian integral $I(h)$ can be expressed as
\begin{eqnarray}
I(h)=\sum\limits_{i=1}^{n}\alpha_ih^i,
\end{eqnarray}
where $\alpha_i,i=1,2,\cdots,n$ are arbitrary constants that can be expressed in terms of $\xi_{i,j}$.}
\vskip 0.2 true cm

\noindent
{\bf Proof}\,  %By inserting (2.17) into $I(h)$  in Lemmas 1 and 2, we immediately obtain equation (2.20). We are now in a position to prove that the coefficients $\alpha_i$ can be chosen as free parameters by mathematical induction.
We prove the proposition by induction on $n$. When $n=2$,  after a direct computation, one has
\begin{eqnarray}I(h)=\xi_{0,1}I_{0,1}(h)+\xi_{0,2}I_{0,2}(h)=\alpha_1h+\alpha_2h^2,\end{eqnarray}
where $$\alpha_1=-\frac{2\pi\sqrt{\lambda}}{\sqrt{1-\lambda}}\xi_{0,1},\ \alpha_2=\frac{2\pi\lambda^\frac32}{(1-\lambda)^\frac32}\xi_{0,2}.$$
By taking advantage of the above two equalities, one gets that the determinant of the following matrix
$$\frac{\partial(\alpha_1,\alpha_2)}{\partial(\xi_{0,1},\xi_{0,2})}=
\begin{pmatrix}
  -\frac{2\pi\sqrt{\lambda}}{\sqrt{1-\lambda}} & 0 \\
  0 & \frac{2\pi\lambda^\frac32}{(1-\lambda)^\frac32} \\
\end{pmatrix}$$
is different from $0$. This implies that $\alpha_1$ and $\alpha_2$ can be taken as free parameters. Similar to the arguments in the case for $n = 2$, when $n=3$, the Abelian integral $I(h)$ can be written as
\begin{eqnarray}\begin{aligned}I(h)&=\xi_{0,1}I_{0,1}(h)+\xi_{0,2}I_{0,2}(h)+\xi_{0,3}I_{0,3}(h)+\xi_{1,2}I_{1,2}(h)+\xi_{2,1}I_{2,1}(h)\\&=\alpha_1h+\alpha_2h^2+\alpha_3h^3,\end{aligned}\end{eqnarray}
where $$\begin{aligned}&\alpha_1=-\frac{2\pi\sqrt{\lambda}}{\sqrt{1-\lambda}}\xi_{0,1},\ \alpha_3=-\frac{3\pi\lambda^\frac52}{(1-\lambda)^\frac52}\xi_{0,3},\\
&\alpha_2=\frac{2\pi\lambda^\frac32}{(1-\lambda)^\frac32}\xi_{0,2}-\frac{3\pi\lambda^\frac32}{(1-\lambda)^\frac32}\xi_{0,3}
+\frac{2\pi\lambda^\frac32}{(1-\lambda)^\frac32}\xi_{1,2}-\frac{\pi\sqrt{\lambda}}{(1-\lambda)^\frac32}\xi_{2,1}.\end{aligned}$$
Note that $\xi_{0,1}$ does not appear in $\alpha_2$, $\xi_{0,1}$ and $\xi_{0,2}$ do not appear in $\alpha_3$ and $\xi_{0,3}$ must appear in $\alpha_3$. Hence, one gets
$$\det\Big[\frac{\partial(\alpha_1,\alpha_2,\alpha_3)}{\partial(\xi_{0,1},\xi_{0,2},\xi_{0,3})}\Big]=
\det\Big[\begin{pmatrix}
  -\frac{2\pi\sqrt{\lambda}}{\sqrt{1-\lambda}} & 0& 0\\
  0 & \frac{2\pi\lambda^\frac32}{(1-\lambda)^\frac32}& -\frac{3\pi\lambda^\frac32}{(1-\lambda)^\frac32}\\
  0 & 0& -\frac{3\pi\lambda^\frac52}{(1-\lambda)^\frac52}
\end{pmatrix}\Big]=\frac{12\pi^3\lambda^\frac92}{(1-\lambda)^\frac92}\neq0,$$
which yields that $\alpha_1$, $\alpha_2$, and $\alpha_3$ can be chosen arbitrarily.
%Now assume the conclusion holds for all $n ¡Ü k-1$. That is,  $\alpha_1$, $\alpha_2,\cdots,\alpha_{k-1}$ can be treated as free parameters. Consequently, the determinant of the following Jacobian matrix
In view of Lemma 2.6 and (2.13), one has
$$\begin{aligned}I(h)=&\sum\limits_{i+j=0}^{n}\xi_{i,j}I_{i,j}(h)\\
=&\sum\limits_{i+j=0}^{n-1}\xi_{i,j}I_{i,j}(h)+\xi_{0,n}I_{0,n}(h)+\xi_{1,n-1}I_{1,n-1}(h)+\cdots+\xi_{n-1,1}I_{n-1,1}(h)\\
=&\tilde{\alpha}_1h+\tilde{\alpha}_{2}h^{2}+\cdots+\tilde{\alpha}_{n-1}h^{n-1}\\&+\xi_{0,n}I_{0,n}(h)+\xi_{1,n-1}I_{1,n-1}(h)+\cdots+\xi_{n-1,1}I_{n-1,1}(h)\\
\triangleq&\sum\limits_{i=1}^{n}\alpha_ih^i,
\end{aligned}$$
where the second and third equalities employ the induction hypothesis. Again, thanks to the induction hypothesis, one has $\tilde{\alpha}_{1}$, $\tilde{\alpha}_{2}$, $\cdots$, $\tilde{\alpha}_{n-1}$ are mutually independent, which implies that the  determinant of the Jacobian matrix
$$\mathbf{A}=\frac{\partial(\alpha_1,\alpha_2,\cdots,\alpha_{n-1})}{\partial(\xi_{i_1,j_1},\xi_{i_2,j_2},\cdots,\xi_{i_{n-1},j_{n-1}})}$$
is non-vanishing, where the sum of the two subscripts of $\xi$ in the above matrix $\mathbf{A}$ is less than $n$. Observe that $\alpha_n$ is the coefficient of $h^n$, hence $\xi_{i_1,j_1},\xi_{i_2,j_2},\cdots,\xi_{i_{n-1},j_{n-1}}$ do not appear in $\alpha_n$. It follows that the partial derivatives of $\alpha_n$ with respect to them all vanish.
According to Lemma 2.6, $\xi_{0,n}$ must appear in $\alpha_n$. Based on the previous arguments, it follows that
$$\mathbf{B}=\frac{\partial(\alpha_1,\alpha_2,\cdots,\alpha_{n-1},\alpha_{n})}{\partial(\xi_{i_1,j_1},\xi_{i_1,j_1},\cdots,\xi_{i_{n-1},j_{n-1}},\xi_{0,n})}
=\left(\begin{array}{ccc}
 \mathbf{A} & \mathbf{C} \\
  \mathbf{0} & \delta  \\
    \end{array}
  \right),$$
where $\mathbf{C}$ is an $(n-1)$-dimensional column vector, $\mathbf{0}$ is an $(n-1)$-dimensional zero row vector, and $\delta$ is a nonzero constant. Therefore, one gets
$$\det(\mathbf{B})=\delta\det(\mathbf{A})\neq0,$$
which implies that $\alpha_{1}$, $\alpha_{2}$, $\cdots$, $\alpha_{n-1},\alpha_n$ are mutually independent. This completes the proof.\quad $\lozenge$
\vskip 0.2 true cm

\noindent
{\bf Remark 2.2}\,  When determining the lower bound of the number of limit cycles, it is essential to verify the independence of the coefficients of $I(h)$ in (2.30). According to conventional methods used in existing literature, this requires obtaining an explicit expression of $\alpha_i$, $i=1,2,\cdots,n$ in terms of $\xi_{i,j}$, $i=0,1,2,\cdots,n;j=0,1,2,\cdots,n$. This is a task that is extremely difficult or even impossible to accomplish. During the proof of Proposition 2.1, by skillfully applying mathematical induction, we successfully circumvent this difficulty by focusing exclusively on the coefficient $\alpha_n$ of $h^n$. This approach not only dramatically reduced computational effort but also achieved what was previously deemed impossible.

\section{Proof of the main result and numerical simulation}
 \setcounter{equation}{0}
\renewcommand\theequation{3.\arabic{equation}}

In order to obtain the lower bound of the number of zeros of $I(h)$, we resort to a result of  Coll, Gasull and Prohens published in \cite{CGP}. We review this result here for the convenience of the reader.
\vskip 0.2 true cm

\noindent
{\bf Lemma 3.1}\, {\it Consider $p+1$ linearly independent analytical functions $f_i:U\rightarrow \mathbb{R},\ i=0,1,2,\cdots, p$, where $U\subset\mathbb{R}$ is an interval. Suppose that there exists $j\in\{0,1,\cdots,p\}$ such that $f_j$ has constant sign. Then there exists $p+1$ constants $\delta_i,\ i=0,1,\cdots,p$,  such that $f(x)=\sum\limits_{i=0}^p\delta_if_i(x)$ has at least $p$ simple zeros in $U$. }
\vskip 0.2 true cm

\noindent
{\bf Proof of Theorem 1.1} When $n=1$, a direct computation gives
$$I(h)=\xi_{0,1}I_{0,1}(h)=-\frac{2\pi\sqrt{\lambda}}{\sqrt{1-\lambda}} h,$$
which yields that system (1.8) has no limit cycle. When $n\geq2$, it follows from Proposition 2.1 that $I(h)$ has at most $n-1$ simple zeros. It is apparent that $h,h^2,\cdots,h^n$ are linearly independent analytical functions and $h$ has constant sign in $(0,+\infty)$. The existence of $n-1$ simple zeros on $(0,+\infty)$ for $I(h)$ can be guaranteed by an appropriate choice of the parameters $\alpha_{1}$, $\alpha_{2}$, $\cdots$, $\alpha_{n-1},\alpha_n$, as indicated by Lemma 3.1. This completes the proof of Theorem 1.1. \quad $\lozenge$
\vskip 0.2 true cm

Next, we will provide corresponding numerical simulations for concrete values of $n$ and $\lambda$ to verify the theoretical result.
 When $n=2$, it follows from Proposition 2.1 and (2.31) that $I(h)$ possesses a simple zero $\frac{(1-\lambda)\xi_{0,1}}{\lambda\xi_{0,2}}$ in $(0,1)$.
Taking $\xi_{0,1}=1,\xi_{0,2}=3$ and $\lambda=\frac{1}{2}$ yields that $I(h)$ has a zero $\frac{1}{3}$. That is, we can find a differential system

When $n=3$, according to (2.32), we take $\alpha_1=\pi,\alpha_2=-5\pi,\alpha_3=6\pi$ and $\lambda=\frac12$, then $I(h)$ has two positive zeros $\frac12$ and $\frac13$. Based on the previous analysis, one can find  a system

When $n=4$, following an analogous approach to the preceding analysis, one has
$$\begin{aligned}
I(h)=\alpha_1h+\alpha_2h^2+\alpha_3 h^3+\alpha_4 h^4,
\end{aligned}$$
where
$$\begin{aligned}&\alpha_1=-\frac{2\pi\sqrt{\lambda}}{\sqrt{1-\lambda}}\xi_{0,1},\ \alpha_4=\frac{5\pi\lambda^\frac72}{(1-\lambda)^\frac72}\xi_{0,4},\\
&\alpha_2=\frac{2\pi\lambda^\frac32}{(1-\lambda)^\frac32}\xi_{0,2}-\frac{3\pi\lambda^\frac32(\lambda-1)}{(1-\lambda)^\frac52}\xi_{0,3}
+\frac{2\pi\lambda^\frac32}{(1-\lambda)^\frac32}\xi_{1,2}-\frac{\pi\sqrt{\lambda}}{(1-\lambda)^\frac32}\xi_{2,1}\\
&\alpha_3=-\frac{3\pi\lambda^\frac52}{(1-\lambda)^\frac52}\xi_{0,3}-\frac{4\pi\lambda^\frac52(2\lambda^2-\lambda-1)}{(1-\lambda)^\frac72}\xi_{0,4}
-\frac{6\pi\lambda^\frac52}{(1-\lambda)^\frac52}\xi_{1,3}+\frac{2\pi\lambda^\frac32}{(1-\lambda)^\frac52}\xi_{2,2}.\end{aligned}$$
Taking $\alpha_1=-4\pi,\alpha_2=55\pi,\alpha_3=-\frac{325}{2}\pi,\alpha_4=125\pi$ and $\lambda=\frac12$, then $I(h)$ has three positive zeros $\frac{1}{10},\frac25$ and $\frac45$. Depending on the values of $\alpha_i(i=1,2,3,4)$, we can identify a system
{\small\begin{eqnarray}\begin{cases}
\frac{dx}{dt}=-2y-x-x^2,\\
\frac{dy}{dt}=x+y+2xy+\frac32x^2+x^3+\varepsilon(2y+\frac{15}{2}xy^2-20x^2y-\frac{245}{8}x^2y^2+40xy^3+ 25y^4)
\end{cases}\end{eqnarray}}
 with three limit cycles, as illustrated in Fig. 4. with $\varepsilon=10^{-4}$.

 \section{Discussion}
 \setcounter{equation}{0}
\renewcommand\theequation{1.\arabic{equation}}

After more than a decade of relentless efforts by numerous scholars, the weak Hilbert's 16th problem has been completely resolved for the case where the Hamiltonian function $H(x,y)$ has degree $\deg H(x,y)=3$ and the perturbation polynomials satisfy $\deg f(x,y)=\deg g(x,y)=2$. This stands as one of the few comprehensive results achieved in this field of research. Bearing this in mind, this paper focuses on the weak Hilbert's 16th problem for a class of cubic isochronous Hamiltonian systems, where the Hamiltonian function is
\begin{eqnarray*}
H(x,y)=\frac{1}{2}x^2+\lambda x^3+\frac{1}{2}\lambda x^4+\frac{1}{2}\lambda^{-1}y^2+xy+x^2y.
\end{eqnarray*}
Under $n$-th degree polynomial perturbations, the exact number of limit cycles is derived.

Unlike previous studies, some terms $I_{i,j}(h)$ appearing in the Abelian integral $I(h)$ of this system cannot be iterated using the derived recurrence formulas. As a result, the number of generators of $I(h)$ depends on the degree $n$ of the perturbation terms. To overcome this difficulty, we identify the differential equations satisfied by these non-iterable  terms and obtain their explicit expressions by solving these differential equations. Another key difficulty lies in verifying the linear independence of coefficients in the expression of $I(h)$ when investigating the lower bound for the number of limit cycles. This obstacle was resolved through mathematical induction. The results presented in this work constitute a meaningful advancement in addressing the weak Hilbert's 16th problem on cubic isochronous Hamiltonian systems. Investigating the weak Hilbert's 16th problem for other types of cubic systems will be an important focus of our future research.
\vskip 0.2 true cm

\noindent
{\bf Acknowledgment}
 \vskip 0.2 true cm

\noindent
This work was supported by the National Natural Science Foundation of China (12161069).
\vskip 0.2 true cm

\noindent
{\bf Author Contributions}
 \vskip 0.2 true cm

\noindent
Jihua Yang: Conceptualization, Funding acquisition, Investigation, Methodology, Project administration, Resources, Supervision, Writing-original draft, Writing-review $\&$ editing.
\vskip 0.2 true cm

\noindent
{\bf Conflict of interest}
 \vskip 0.2 true cm

\noindent
The authors declare that they have no known competing financial interests or personal relationships that could have appeared to influence the work reported in this paper.
\vskip 0.2 true cm

\noindent
{\bf Data availability statement}
 \vskip 0.2 true cm

\noindent
No data was used for the research in this article. It is pure mathematics.

\end{document}